\documentclass[a4, 12pt]{amsart}
\usepackage{amsmath}
\usepackage{amssymb}
\numberwithin{equation}{section}
\usepackage{graphicx}
\newtheorem{thm}{Theorem}[section]

\newtheorem{example}[thm]{Example}
\newtheorem{theorem}{Theorem}[section]

\newtheorem{remark}{Remark}[section]
\numberwithin{equation}{section}

\begin{document}

\title[Estimates for the first eigenvalue of Jacobi operator]
{Estimates for the first eigenvalue of Jacobi operator on  hypersurfaces with constant mean curvature in  spheres}
\author{Daguang Chen and Qing-Ming Cheng}
\address{Daguang Chen\\
Department of Mathematical Sciences,
Tsinghua University,
Beijing 100084, P. R. China, dgchen@math.tsinghua.edu.cn
}
\address{Qing-Ming Cheng \\ Department of Applied Mathematics, Faculty of Sciences,
Fukuoka  University, 814-0180, Fukuoka,  Japan, cheng@fukuoka-u.ac.jp}
\thanks{The first author is supported by NSFC:  No.11471180.
The second author is partially  supported by JSPS Grant-in-Aid for Scientific Research (B): No.16H03937.}

\date{}%
\begin{abstract}
In this paper, we study the first eigenvalue of Jacobi operator on
an $n$-dimensional non-totally umbilical compact hypersurface with
constant mean curvature $H$  in the  unit sphere $S^{n+1}(1)$. 
We give an optimal upper bound for  the first eigenvalue of Jacobi operator,
which only depends on the mean curvature $H$ and the dimension $n$. 
This bound is attained  if and only if,  $\varphi:\  M \to S^{n+1}(1)$ is isometric to 
$S^1(r)\times S^{n-1}(\sqrt{1-r^2})$ when $H\neq 0$ or  $\varphi:\  M \to S^{n+1}(1)$ is isometric to 
a Clifford torus $
S^{n-k}(\sqrt{\dfrac{n-k}n})\times S^k(\sqrt{\dfrac{k}n})$,
 for $k=1,  2, \cdots, n-1$ when $H=0$.
 \end{abstract}
\maketitle

\section{Introduction}

\noindent
Let $\varphi:  M \to S^{n+1}(1)$ be an $n$-dimensional compact hypersurface in the unit sphere $S^{n+1}(1)$ 
of dimension  $n+1$.  We consider
 a variation of the hypersurface  $\varphi:  M \to S^{n+1}(1)$, for any  $t\in (-\varepsilon, \varepsilon)$, 
$$
\varphi_t:  M \to S^{n+1}(1)
$$
is an immersion with
$\varphi_0=\varphi$. 
The area of $\varphi_t$ is given by 
$$
A(t)=\int_M dA_t
$$
and  the volume of $\varphi_t$ is defined by 
$$
V(t)=\dfrac1{n+1}\int_{M}\langle \varphi_t, N(t)\rangle dA_t,
$$
where $N(t)$ denotes the unit normal of $\varphi_t$.
For any $t$, if $V(t)=V(0)$, then the variation 
$\varphi_t$ is called volume-preserving. If the variational
vector $\dfrac{\partial \varphi_t}{\partial t}|_{t=0}=f N$ for a
smooth function $f$, then the variation is called a
normal variation, where $N$ is the unit normal of 
$\varphi$.
Let $H$ denote the mean curvature of $\varphi$. The first variation
formula of the area functional $A(t)$ is given by
$$
\frac{dA(t)}{dt}|_{t=0}=-\int_MnHfdA,
$$
where $f=\langle \frac{\partial \varphi_t}{\partial t}|_{t=0}, N\rangle.$
Thus, we know that a compact  hypersurface is minimal, that is, $H\equiv 0$
if and only if 
$$
\frac{dA(t)}{dt}|_{t=0}=0.
$$
Hence, compact minimal hypersurfaces are critical points of the area functional
$A(t)$.
The second variation formula of  $A(t)$ is given by 
$$
\dfrac{d^2A(t)}{dt^2}|_{t=0}=-\int_MfJfdA
$$
and 
$$
Jf=\Delta f+(S+n)f,
$$
where $S$ denotes the squared norm of the second fundamental form
of  $\varphi$ 
and  $\Delta$ stands for the Laplace-Beltrami operator.
The $J$ is called a Jacobi operator or  a stability operator on the minimal hypersurface $\varphi$ (cf. \cite{a}, \cite{p}). 

\noindent
Let $\lambda^{J}_1$ denote the first eigenvalue of the Jacobi operator $J$. Then
$$
Ju=-\lambda^{J}_1 u
$$
and 
the $\lambda^{J}_1$ is given by
$$
\lambda^{J}_1=\inf _{f\not\equiv 0}\dfrac{-\int_MfJfdA}{\int_Mf^2dA}.
$$
For a compact minimal hypersurface in $S^{n+1}(1)$,  Simons \cite{s} proved 
$$
\lambda^{J}_1\leq -n
$$
and  $\lambda^{J}_1=-n$ if and only if  $\varphi:  M \to S^{n+1}(1)$  is
totally geodesic.
Furthermore, 
Wu \cite {w}  proved that for an $n$-dimensional compact non-totally geodesic 
minimal hypersurface   $\varphi:  M \to S^{n+1}(1)$  in $S^{n+1}(1)$,
then  $\lambda^{J}_1 \leq -2n$ and $\lambda^{J}_1=-2n$ if and only if  $\varphi:  M \to S^{n+1}(1)$  is
a Clifford torus
$
S^{n-k}(\sqrt{\frac{n-k}n})\times S^k(\sqrt{\frac{k}n}),
$
 for $k=1, 2, \cdots, n-1$. Thus, we know that the upper bound for the first eigenvalue  
 $\lambda^{J}_1$ due to Wu is optimal and it only depends on the dimension $n$, does not depends on the immersion. 

\vskip3mm
\noindent
On the other hand, if one  considers  the volume-preserving variation of $\varphi$,
then we have 
$$
\int_MfdA=0.
$$
From the first variation formula:
$$
\dfrac{dA(t)}{dt}|_{t=0}=-\int_MnHfdA,
$$
we know that  compact hypersurfaces with constant mean curvature 
are critical points of the area functional $A(t)$ for  the volume-preserving variation
 and the second variation formula of  $A(t)$ is given by
$$
\dfrac{d^2A(t)}{dt^2}|_{t=0}=-\int_MfJfdA,
$$
where the Jacobi operator $J$ on compact hypersurfaces with constant mean curvature 
is the same as one of compact minimal hypersurfaces (\cite{a}, \cite{bce}).

\vskip 2mm
\noindent
Alias, Barros and Brasil \cite{abb} studied the first eigenvalue of the Jacobi operator
$J$ on compact hypersurfaces with constant mean curvature. They proved the following:

\vskip3mm
\noindent
{\bf Theorem ABB}. {\it 
If  $\varphi:  M \to S^{n+1}(1)$ is an $n$-dimensional compact 
hypersurface with non-zero constant mean curvature $H$ in the unit sphere 
$S^{n+1}(1)$, then either 
$\lambda^{J}_1=-n(1+H^2)$ and $\varphi:   M \to S^{n+1}(1)$ is totally umbilical or
\begin{equation*}
\begin{aligned}
\lambda^{J}_1\leq &-2n(1+H^2)+\dfrac{n(n-2)|H|}{\sqrt{n(n-1)}}\max \sqrt{S-nH^2}
\end{aligned}
\end{equation*}
and  the equality holds if and only if $\varphi:   M \to S^{n+1}(1)$ is 
$S^1(r)\times S^{n-1}(\sqrt{1-r^2})$, with $r^2>\frac{1}n$ for $n\geq 2$.
}

\vskip2mm

\noindent
According to this theorem, we know that,  for $n=2$, the upper bund of the first eigenvalue
$\lambda^{J}_1$ of the Jacobi operator  of non-totally umbilical  compact hypersurfaces with constant mean curvature
only depends on the mean curvature $H$ and the dimension. But for $n\geq 3$, 
the upper bound of the first eigenvalue
$\lambda^{J}_1$ of the Jacobi operator  on non-totally umbilical  compact hypersurfaces with constant mean curvature
includes  the term $\max \sqrt{S-nH^2}$. Hence,
the upper bound of the first eigenvalue
$\lambda^{J}_1$  does not only depend  on the mean curvature $H$ and the dimension $n$, but also depends on the immersion $\varphi$. 

\noindent
It is natural and important to propose the following:

\vskip2mm
\noindent
{\bf Problem 1.1}.  To find an optimal upper bound for  the first eigenvalue
$\lambda^{J}_1$ of the Jacobi operator  on non-totally umbilical  compact 
hypersurfaces with constant mean curvature, which 
only depends on the mean curvature $H$ and the dimension $n$.

\vskip2mm
\noindent
In this paper, we give an affirmative answer for the above problem 1.1.

\begin{theorem}
Let $\varphi: M \to S^{n+1}(1)$ be  an $n$-dimensional non-totally umbilical compact 
hypersurface with constant mean curvature $H$ in the unit sphere $S^{n+1}(1)$.
\begin{enumerate}
\item If  $2\leq n\leq 4 $ or $n\geq 5$ and $n^2H^2< \frac{16(n-1)}{n(n-4)}$, 
then the first eigenvalue $\lambda_1^{J}$ of  the Jacobi operator $J$ satisfies 
$$
\lambda_1^{J}\leq -n(1+H^2)-\dfrac{n(\sqrt{4(n-1)+n^2H^2}-(n-2)|H|)^2}{4(n-1)}
$$
and  the equality holds if and only if  $\varphi: M \to S^{n+1}(1)$ is isometric to 
$S^1(r)\times S^{n-1}(\sqrt{1-r^2})$ with  $r>0$ satisfying
$$
\begin{cases}
1>r^2>\dfrac{1}n & \ \text{for}  \ 2\leq n\leq 4,\\ 
 \dfrac{n}{(n-2)^2}>r^2>\dfrac{1}n, &\ \text{for} \ n\geq 5  \ \text{and}  \  n^2H^2< \dfrac{16(n-1)}{n(n-4)}
\end{cases}
$$
or  $\varphi:  M \to S^{n+1}(1)$ is isometric to 
a Clifford torus
$
S^{n-k}(\sqrt{\frac{n-k}n})\times S^k(\sqrt{\frac{k}n}),
$
 for $k=1, 2, \cdots, n-1$ with $H=0$.
\item If  $n\geq 5$ and $n^2H^2\geq\frac{16(n-1)}{n(n-4)}$, the first eigenvalue $\lambda_1^{J}$ of  the Jacobi operator $J$ satisfies
\begin{equation*}
\begin{aligned}
\lambda^{J}_1\leq -2(n-1)(1+H^2)+\dfrac{(n-2)^4}{8(n-1)}H^2
\end{aligned}
\end{equation*}
and  the equality holds if and only if  $\varphi: M \to S^{n+1}(1)$ is isometric to 
$S^1(\frac{\sqrt n}{n-2})\times S^{n-1}(\frac{\sqrt{(n-1)(n-4)}}{n-2})$.
\end{enumerate}
\end{theorem}

\begin{remark} 
Since  the first eigenvalue  of Jacobi operator $J$  on 
totally umbilical hypersurfaces satisfies $\lambda_1^{J}=-n(1+H^2)$, according to our theorem, one  knows 
that for $2\leq n\leq 4$, there are no $n$-dimensional compact hypersurfaces in the unit sphere with constant mean
curvature $H$ so that the first eigenvalue $\lambda_1^{J}$ of Jacobi operator $J$ takes a value in the internal
$$
\biggl(-n(1+H^2)-\dfrac{n(\sqrt{4(n-1)+n^2H^2}-(n-2)|H|)^2}{4(n-1)}, \ -n(1+H^2)\biggl).
$$
For any $n\geq 2$, there are no $n$-dimensional compact hypersurfaces in the unit sphere with constant mean
curvature $H$ satisfying $n^2H^2<\frac{16(n-1)}{n(n-4)}$
so that the first eigenvalue $\lambda_1^{J}$ of Jacobi operator $J$ takes a value in the internal
$$
\biggl(-n(1+H^2)-\dfrac{n(\sqrt{4(n-1)+n^2H^2}-(n-2)|H|)^2}{4(n-1)}, \ -n(1+H^2)\biggl).
$$
One should compare the bound  
$$
 -n(1+H^2)-\dfrac{n(\sqrt{4(n-1)+n^2H^2}-(n-2)|H|)^2}{4(n-1)}
 $$
with the pinching constant in the rigidity theorem of Cheng and Nakagawa \cite{cn} or  Alencar and  do Carmo \cite{ac}.  
\end{remark}

\vskip3mm
\noindent
{\bf Acknowledgement.}   The authors would like to express their thanks to the referee for the valuable comments and suggestions.

\section{Preliminaries}

\vskip3mm
\noindent
Throughout this paper, all manifolds are assumed to be smooth and
connected without boundary. Let $\varphi:   M \to S^{n+1}(1)$ be an
$n$-dimensional hypersurface
in a unit sphere $S^{n+1}(1)$. We choose a local orthonormal frame
$\{{\bf e}_1, \cdots, {\bf e}_{n}, {\bf e}_{n+1}\}$ and the dual
coframe $\{\omega_1, \cdots,$  $\omega_n$,  $ \omega_{n+1}\}$ in
such a way that $\{{\bf e}_1, \cdots, {\bf e}_n\}$  is a local
orthonormal frame on $M$. Hence, we have
$$
\omega_{n+1}=0
$$
on $M$. From Cartan's lemma, we have
\begin{equation}
\omega_{i n+1}=\sum_{j=1}^nh_{ij}\omega_j, \ h_{ij}=h_{ji}. 
\end{equation}
The mean curvature $H$ and the second fundamental form
${II}$ of $\varphi:\  M \to S^{n+1}(1)$ are defined, respectively, by
$$
H=\frac{1}{n}\sum_{i=1}^nh_{ii}, \
{II}=\sum_{i,j=1}^nh_{ij}\omega_i\otimes\omega_j{\bf e}_{n+1}.
$$
When the mean curvature $H$ of $\varphi:  M \to S^{n+1}(1)$ is identically zero, we
recall that $\varphi:   M \to S^{n+1}(1)$ is by definition {\it a minimal hypersurface}.
From the structure equations of $\varphi:   M \to S^{n+1}(1)$, Gauss equation is  given by
\begin{equation}
R_{ijkl}=(\delta_{ik}\delta_{jl}-\delta_{il}\delta_{jk})
+(h_{ik}h_{jl}-h_{il}h_{jk}), 
\end{equation}
 From
(2.2), we have
\begin{align*}
 n(n-1)r=n(n-1)+n^2H^2-S,
\end{align*}
where $n(n-1)r$ and $S$ denote the scalar curvature and the squared norm
of the second fundamental form of $\varphi:   M \to S^{n+1}(1)$, respectively. 
Defining the covariant derivative of $h_{ij}$ by
\begin{equation}
\sum_{k}h_{ijk}\omega_k=dh_{ij}+\sum_k
h_{ik}\omega_{kj} +\sum_k h_{kj}\omega_{ki},
\end{equation}
we obtain the Codazzi equations

\begin{equation}
h_{ijk}=h_{ikj}.
\end{equation}
By taking exterior differentiation of (2.3), and defining
\begin{equation}
\sum_lh_{ijkl}\omega_l=dh_{ijk}+\sum_lh_{ljk}\omega_{li}
+\sum_lh_{ilk}\omega_{lj}+\sum_l
h_{ijl}\omega_{lk},
\end{equation}
we have the following Ricci identities:
\begin{equation}
h_{ijkl}-h_{ijlk}=\sum_m
h_{mj}R_{mikl}+\sum_m h_{im}R_{mjkl}.
\end{equation}

\noindent
For any $C^2$-function $f$ on $M$, we define  its gradient and Hessian by
$$
df=\sum_{i=1}^nf_i\omega_i,
$$
$$
\sum_{j=1}^nf_{ij}\omega_j=df_i+\sum_{j=1}^nf_j\omega_{ji}.
$$
Thus, the Laplace-Beltrami operator  $\Delta $ is given by
$$
\Delta f=\sum_{i=1}^nf_{ii}.
$$

\begin{example}
For totally umbilical sphere $S^n(r)$ of radius $r>0$, the first eigenvalue $\lambda_1^{J}=-n(1+H^2)$ with $H=\frac1r$.
\end{example}
\begin{example}
For Clifford torus $S^{n-k}(\sqrt{\frac{n-k}n})\times S^k(\sqrt{\frac{k}n})$, $k=1, 2, \dots, n$, 
the first eigenvalue $\lambda_1^{J}=-2n$ with $H=0$.
\end{example}
\begin{example}
For hypersurfaces
$S^1(r)\times S^{n-1}(\sqrt{1-r^2})$ with $0<r<1$,  the principal curvatures are given by
$$
k_1=-\dfrac{\sqrt{1-r^2}}r, \  \ k_2=\cdots=k_n=\dfrac{r}{\sqrt{1-r^2}}.
$$
Hence, we know that 
$$
nH=\dfrac{nr^2-1}{r\sqrt{1-r^2}}, \ \ S=\dfrac{1-2r^2+nr^4}{r^2(1-r^2)}.
$$
For  $r^2\geq \frac1n$, by a direct computation, we know that the first eigenvalue 
$\lambda_1^{J}$ of the Jacobi operator $J$ on $S^1(r)\times S^{n-1}(\sqrt{1-r^2})$ satisfies
$$
 \lambda_1^{J}=-n(1+H^2)-\dfrac{n(\sqrt{4(n-1)+n^2H^2}-(n-2)|H|)^2}{4(n-1)}.
$$
For $n\ge 5$ and   $\frac{1}n\leq r^2< \frac{n}{(n-2)^2}$, we know the hypersurface $S^1(r)\times S^{n-1}(\sqrt{1-r^2})$
satisfies 
$$
n^2H^2< \frac{16(n-1)}{n(n-4)}
$$
and 
$$
 \lambda_1^{J}=-n(1+H^2)-\dfrac{n(\sqrt{4(n-1)+n^2H^2}-(n-2)|H|)^2}{4(n-1)}.
$$
The hypersurface $S^1(\frac{\sqrt n}{n-2})\times S^{n-1}(\frac{\sqrt{(n-1)(n-4)}}{n-2})$ satisfies

\begin{equation*}
\begin{aligned}
\lambda^{J}_1=-2(n-1)(1+H^2)+\dfrac{(n-2)^4}{8(n-1)}H^2
\end{aligned}
\end{equation*}
with  $n^2H^2=\frac{16(n-1)}{n(n-4)}$.

\end{example}

\section{Proof of  theorem 1.1.}

\vskip3mm
\noindent
In this section, we give a proof of the theorem 1.1.

\vskip2mm
\noindent
{\it Proof of theorem} 1.1.
When  $H\equiv 0$, according to the result of Wu \cite{w},
we have  $\lambda_1^{J}\leq -2n$ and $\lambda_1^{J}= -2n$ if and only if 
 $\varphi: M \to S^{n+1}(1)$ is isometric to 
a Clifford torus 
$S^{n-k}(\sqrt{\frac{n-k}n})\times S^k(\sqrt{\frac{k}n})$,
 for $k=1,  2, \cdots, n-1$.

\noindent
From now we assume $H\neq 0$.
By making use of the Codazzi equations,  Ricci identities  and a standard computation of Simons'  type formula
(cf. \cite{cn},  \cite{c1, c2}, \cite{lh} and \cite{s}),
we have 
\begin{equation}
\frac12\Delta S=\sum_{i,j,k=1}^nh_{ijk}^2+nS-n^2H^2+nHf_3-S^2,
\end{equation}
where $f_3=\sum_{i=1}^nk_i^3$ and $k_i$, $i=1, 2, \dots, n$ denote the principal curvatures.

\noindent
Putting
$\mu_i=k_i-H$, we have 
\begin{equation}
B:=\sum_{i=1}^n\mu_i^2=S-nH^2\geq 0, \  \ f_3=B_3+3HB +nH^3,
\end{equation}
where $B_3=\sum_{i=1}^n\mu_i^3$.
The following inequality is known (cf. \cite{cn} and \cite{lh}):
\begin{equation}
|B_3|\leq \dfrac{n-2}{\sqrt{n(n-1)}}B^{\frac32},
\end{equation}
and the equality holds if and only if  at least $n-1$  of $k_i$, for $i=1, 2, \dots, n$, are equal with each other.
Since $H$ is constant, we can assume $H> 0$. Thus, from (3.1), (3.2) and (3.3), we have
\begin{equation}
\begin{aligned}
\frac12\Delta B=\frac12\Delta S\geq \sum_{i,j,k=1}^nh_{ijk}^2+B(n+nH^2-B)-nH\dfrac{n-2}{\sqrt{n(n-1)}}B^{\frac32}.
\end{aligned}
\end{equation}
For any constant  $\alpha>0$ and $\varepsilon>0$, we consider a function $f_{\varepsilon}=(B+\varepsilon)^{\alpha}>0$.
Hence, we have, from (3.4),  
\begin{equation}
\begin{aligned}
\Delta f_{\varepsilon}&=\alpha(\alpha-1)(B+\varepsilon)^{\alpha-2}|\nabla B|^2+\alpha (B+\varepsilon)^{\alpha-1}\Delta B\\
&\geq \alpha(\alpha-1)(B+\varepsilon)^{\alpha-2}|\nabla B|^2\\
&+2\alpha (B+\varepsilon)^{\alpha-1}\biggl(\sum_{i,j,k=1}^nh_{ijk}^2+B(n+nH^2-B)-nH\dfrac{n-2}{\sqrt{n(n-1)}}B^{\frac32}\biggl).
\end{aligned}
\end{equation}
Since $H$ is constant, we have
\begin{equation}
\begin{aligned}
&\nabla_k(nH)=\sum_{i=1}^nh_{iik}=0, \  \  \  h_{kkk}^2\leq (n-1)\sum_{i\neq k}h_{iik}^2\\
&|\nabla B|^2=\sum_{k=1}^n(2\sum_{i=1}^n\mu_ih_{iik})^2\leq 4B\sum_{i,k=1}^nh_{iik}^2.
\end{aligned}
\end{equation}
Thus, we obtain
\begin{equation}
\begin{aligned}
|\nabla B|^2&\leq 4B\sum_{i,k=1}^nh_{iik}^2\\
&=4B\bigl(\dfrac{n}{n+2}\sum_{k=1}^nh_{kkk}^2+\dfrac{2}{n+2}\sum_{k=1}^nh_{kkk}^2+\sum_{i\neq k}h_{iik}^2\bigl)\\
&\leq \dfrac{4n}{n+2}B\bigl(\sum_{k=1}^nh_{kkk}^2+3\sum_{i\neq k}h_{iik}^2\bigl).
\end{aligned}
\end{equation}
For any constant $\beta$, we have 
\begin{equation*}
\begin{aligned}
&\lambda_1^{J} \int_Mf_{\varepsilon}^2dA
\leq -\int_Mf_{\varepsilon}J f_{\varepsilon}dA\\
\end{aligned}
\end{equation*}
\begin{equation*}
\begin{aligned}&=-\beta\int_{M}f_{\varepsilon}\Delta f_{\varepsilon}dA
-\int_{M}\biggl((1-\beta)f_{\varepsilon}\Delta f_{\varepsilon}+(S+n)f_{\varepsilon}^2\biggl)dA\\
&=\beta\int_{M}|\nabla f_{\varepsilon}|^2dA-\int_{M}f_{\varepsilon}\biggl\{(1-\beta)\biggl(\alpha(\alpha-1)(B+\varepsilon)^{\alpha-2}|\nabla B|^2\\
&+\alpha (B+\varepsilon)^{\alpha-1}\Delta B\biggl)+(B+nH^2+n)f_{\varepsilon}\biggl\}dA\\
&=\alpha \int_{M}f_{\varepsilon}\bigl\{1+2\alpha\beta-\beta-\alpha\bigl\}(B+\varepsilon)^{\alpha-2}|\nabla B|^2dA\\
&-\int_{M}f_{\varepsilon}^2\biggl\{\dfrac{\alpha(1-\beta)}{B+\varepsilon}
\Delta B+B+nH^2+n\biggl\}dA.
\end{aligned}
\end{equation*}
By  taking  $\alpha$ and $\beta$  satisfying
\begin{equation}
\alpha >\frac{n-2}{4n}, \  \  1-\beta= \dfrac{2n\alpha}{4n\alpha +2-n},
\end{equation}
we have 
$$
(n-2)(1-\beta)-4n\alpha(1-\beta)
+2n\alpha=0.
$$
Since
$$
\sum_{i,j,k=1}^nh_{ijk}^2=\sum_{k=1}^nh_{kkk}^2+3\sum_{i\neq k}h_{iik}^2+\sum_{i\neq j\neq k\neq i}^nh_{ijk}^2,
$$
from (3.7), 
we obtain
\begin{equation}
\begin{aligned}
&\bigl(1+2\alpha\beta-\beta-\alpha\bigl)|\nabla B|^2
-2(1-\beta) (B+\varepsilon)\sum_{i,j,k=1}^nh_{ijk}^2\\
&\leq \dfrac{2}{n+2} B\biggl\{(n-2)(1-\beta)-4n\alpha(1-\beta)
+2n\alpha \biggl\}\bigl(\sum_{k=1}^nh_{kkk}^2+3\sum_{i\neq k}h_{iik}^2\bigl)=0.
\end{aligned}
\end{equation}
Thus,
we infer
\begin{equation*}
\begin{aligned}
&\lambda_1^{J} \int_Mf_{\varepsilon}^2dA\\
&\leq\alpha \int_{M}f_{\varepsilon}(B+\varepsilon)^{\alpha-2}\biggl\{\bigl(1
+2\alpha\beta-\beta-\alpha\bigl)|\nabla B|^2-2(1-\beta) (B+\varepsilon)\sum_{i,j,k=1}^nh_{ijk}^2\biggl\}dA\\
&-\int_{M}f_{\varepsilon}^2\biggl\{\dfrac{2\alpha(1-\beta)B}{B+\varepsilon }
\biggl((n+nH^2-B)- nH\dfrac{(n-2)}{\sqrt{n(n-1)}}B^{\frac12}\biggl)+B+nH^2+n\biggl\}dA\\
&\leq-\int_{M}f_{\varepsilon}^2\dfrac{B}{B+\varepsilon}
\biggl(\bigl\{1-2\alpha(1-\beta)\bigl\}B
-\dfrac{2\alpha(1-\beta)(n-2)}{\sqrt{n(n-1)}} nHB^{\frac12}+\varepsilon \biggl)dA\\
&-2\alpha(1-\beta)(n+nH^2)\int_{M}f_{\varepsilon}^2\dfrac{B}{B+\varepsilon}dA
-(n+nH^2)\int_{M}f_{\varepsilon}^2dA.\\
\end{aligned}
\end{equation*}
For $1-2\alpha(1-\beta)> 0$, we obtain
\begin{equation*}
\begin{aligned}
&\lambda_1^{J} \int_Mf_{\varepsilon}^2dA\\
&\leq\int_{M}f_{\varepsilon}^2\dfrac{B}{B+\varepsilon}
\biggl(\dfrac{\alpha^2(1-\beta)^2(n-2)^2}{(1-2\alpha(1-\beta))n(n-1)} (nH)^2-\varepsilon \biggl)dA\\
&-2\alpha(1-\beta)(n+nH^2)\int_{M}f_{\varepsilon}^2\dfrac{B}{B+\varepsilon}dA
-(n+nH^2)\int_{M}f_{\varepsilon}^2dA.\\\end{aligned}
\end{equation*}
Since $\varphi: M \to S^{n+1}(1)$ is not totally umbilical, we have 
$$
\lim_{\varepsilon\to 0}\int_M f_{\varepsilon}^2dA=\int_MB^{2\alpha}dA>0.
$$
Letting $\varepsilon\to 0$, 
we derive
\begin{equation}
\begin{aligned}
&\lambda^{J}_1 \leq -(1+2\alpha(1-\beta))n(1+H^2)
+\dfrac{\alpha^2(1-\beta)^2}{\big(1-2\alpha(1-\beta)\big)}\dfrac{(n-2)^2}{n(n-1)}n^2H^2.
\end{aligned}
\end{equation}
For $n=2$, 
we have 
\begin{equation*}
\begin{aligned}
&\lambda^{J}_1 \leq -(1+2\alpha(1-\beta))n(1+H^2).
\end{aligned}
\end{equation*} 
From (3.8), we have $\beta =\dfrac12$ for any  $0<\alpha< 1$. Hence, we obtain
\begin{equation*}
\begin{aligned}
&\lambda^{J}_1 \leq -2n(1+H^2).
\end{aligned}
\end{equation*}
For $2< n\leq 4$  or  $n\geq 5$ and $n^2H^2<\dfrac{16(n-1)}{n(n-4)}$, 
we have 
\begin{equation}
\frac12>\dfrac12\biggl(1-\sqrt{\dfrac{(n-2)^2H^2}{4(n-1)+n^2H^2}}\biggl)>  \dfrac12-\dfrac1n\geq  \dfrac12-\dfrac1{\sqrt{2n}}.
\end{equation}
Observe from (3.8) that $1-2\alpha (1-\beta)>0$ if and only if 
\begin{equation}
\frac12-\frac1{\sqrt{2n}}<\alpha<\frac12+\frac1{\sqrt{2n}}.
\end{equation}
Defining
$$
w(\alpha)=\alpha(1-\beta)=\dfrac{2n\alpha^2}{4n\alpha +2-n},
$$
$w(\alpha)$ is an increasing function of $\alpha$, for $\alpha> \frac12-\frac1n$
and
$$
w( \frac12-\frac1n)=\frac12-\frac1n, \quad w(\dfrac12+\dfrac1{\sqrt{2n}})=\frac12.
$$
According to (3.11) and  (3.12),  there exists a $\alpha$ satisfying 
$$
\frac12-\frac1{n}<\alpha
<\frac12+\frac1{\sqrt{2n}}
$$
such that 
\begin{equation}
w(\alpha)= \dfrac12\biggl(1-\sqrt{\dfrac{(n-2)^2H^2}{4(n-1)+n^2H^2}}\biggl).
\end{equation}
Therefore, we have, for this $\alpha$,  
\begin{equation}
1-2\alpha(1-\beta)= \sqrt{\dfrac{(n-2)^2H^2}{4(n-1)+n^2H^2}}>0.
\end{equation}
From (3.10),
we obtain
\begin{equation}
\begin{aligned}
&\lambda^{J}_1\leq  -n(1+H^2)\\
&\quad -2\alpha(1-\beta)n\dfrac{4(n-1)\big(1-2\alpha(1-\beta)\big)(1+H^2)
-2\alpha(1-\beta)(n-2)^2H^2}{4(n-1)\big(1-2\alpha(1-\beta)\big)}.
\end{aligned}
\end{equation}
From (3.14), we infer
\begin{equation*}
\begin{aligned}
&4(n-1)\big(1-2\alpha(1-\beta)\big)(1+H^2)
-2\alpha(1-\beta)(n-2)^2H^2\\
&=\biggl\{4(n-1)(1+H^2)\sqrt{\dfrac{(n-2)^2H^2}{4(n-1)+n^2H^2}} -
\biggl(1-\sqrt{\dfrac{(n-2)^2H^2}{4(n-1)+n^2H^2}}\biggl)(n-2)^2H^2\biggl\}\\
&=4(n-1)\sqrt{\dfrac{(n-2)^2H^2}{4(n-1)+n^2H^2}}
 - (n-2)^2H^2+
\sqrt{\dfrac{(n-2)^2H^2}{4(n-1)+n^2H^2}}n^2H^2\\
&=\sqrt{4(n-1)+n^2H^2}\sqrt{(n-2)^2H^2}
- (n-2)^2H^2\\
&=\sqrt{(n-2)^2H^2}\biggl(\sqrt{4(n-1)+n^2H^2}
- (\sqrt{(n-2)^2H^2}\biggl).\\
\end{aligned}
\end{equation*}
From (3.14), (3.15) and the above equality, we obtain
\begin{equation*}
\begin{aligned}
&\lambda^{J}_1\leq  
-n(1+H^2)-\dfrac{n\biggl(1-\sqrt{\dfrac{(n-2)^2H^2}{4(n-1)+n^2H^2}}\biggl)}{
4(n-1)\sqrt{\dfrac{(n-2)^2H^2}{4(n-1)+n^2H^2}}}\times\\
& \quad \times \sqrt{(n-2)^2H^2}\biggl(\sqrt{4(n-1)+n^2H^2}
- (\sqrt{(n-2)^2H^2}\biggl)\\
&= -n(1+H^2)-\dfrac{n}{4(n-1)}(\sqrt{4(n-1)+n^2H^2}-(n-2)|H|)^2.
\end{aligned}
\end{equation*}
 If  the equality holds, 
we know that  $h_{ijk}=0$, for any $i, j, k=1, 2, \dots, n$. Hence, we know that 
the second fundamental form is parallel and  $S$ is constant. Thus, we know that 
$\varphi: M \to S^{n+1}(1)$ is isometric to 
 $S^1(r)\times S^{n-1}(\sqrt{1-r^2})$ since, from the (3.3),  the $n-1$ of the principal curvatures are equal with each other.
From the examples in the section 2, we know that $r$ satisfies
$$
\begin{cases}
r^2>\dfrac{1}n & \ \text{for}  \ 2\leq n\leq 4,\\ 
\dfrac{1}n< r^2< \dfrac{n}{(n-2)^2}, &\ \text{for} \ n\geq 5  \ \text{and}  \  n^2H^2< \dfrac{16(n-1)}{n(n-4)}.
\end{cases}
$$
\vskip2mm
\noindent
If $n\geq 5$ and $n^2H^2\geq\frac{16(n-1)}{n(n-4)}$, we take 
$$
\alpha(1-\beta)= \dfrac12-\dfrac1n,
$$
that is,
$$
\beta=0 \quad \text{\rm and}\quad  \alpha= \dfrac12-\dfrac1n,
$$ 
Thus, the inequality (3.10) becomes 
$$
\lambda^{J}_1\leq -2(n-1)(1+H^2)+\dfrac{(n-2)^4}{8(n-1)}H^2.
$$
If the equality holds, 
we know 
$$
(1-2\alpha)\sqrt B
=\dfrac{\alpha(n-2)}{\sqrt{n(n-1)}}nH.
$$
Thus, we have 
\begin{equation}
S=B+nH^2=nH^2+\dfrac{(n-2)^4}{16n(n-1)}n^2H^2.
\end{equation}
because of 
$$
\alpha= \dfrac12-\dfrac1n.
$$
Since $S$ is constant,  the first eigenvalue $\lambda_1^{J}$ of the Jacobi operator  is given
by
$$
\lambda_1^{J}=-S-n= -2(n-1)(1+H^2)+\dfrac{(n-2)^4}{8(n-1)}H^2.
$$
Hence, we obtain
\begin{equation}
S= n-2+ 2(n-1)H^2-\dfrac{(n-2)^4}{8(n-1)}H^2.
\end{equation}
From (3.16) and (3.17), we get
\begin{equation*}
n-2=(2-n)H^2+\dfrac{(n-2)^4(n+2)}{16(n-1)}H^2,
\end{equation*}
\begin{equation*}
1=\dfrac{n(n-4)}{16(n-1)}n^2H^2,
\end{equation*}
that is,
\begin{equation*}
n^2H^2=\dfrac{16(n-1)}{n(n-4)}.
\end{equation*}
Since, from the (3.3),  the $n-1$ of the principal curvatures are equal with each other,
From the examples in the section 2,  
we know  that $\varphi: M \to S^{n+1}(1)$ is isometric to 
$S^1(\frac{\sqrt n}{n-2})\times S^{n-1}(\frac{\sqrt{(n-1)(n-4)}}{n-2})$. 
It completes the proof of theorem 1.1.
\begin{flushright}
$\square$
\end{flushright}
\vskip 5mm

\vskip5mm

\bibliographystyle{amsplain}

\end{document}